\documentclass[12pt]{article}
\usepackage{amsmath,amssymb,amsthm}
\usepackage{fullpage}

\newtheorem{thm}{\textbf{Theorem}}

\newcommand{\dis}{\displaystyle}

\begin{document}

\title{The Probability of a Run }
\author{Mark B. Villarino\\
Depto.\ de Matem\'atica, Universidad de Costa Rica,\\
2060 San Jos\'e, Costa Rica}
\date{\today}

\maketitle

 \begin{abstract}
 We deduce the explicit formula for the probability of a run of $r$ successes in $n$ trials.
\end{abstract}
\section{Introduction}

A famous problem in classical probability was first stated in \textsc{De Moivre}'s treatise, \emph{The Doctrine of Chances} \cite{D} as \emph{Problem LXXIV}:\begin{quotation}``To find the Probability of throwing a Chance assigned a given number of times without intermission, in any given number of Trials."\end{quotation}
We formulate this more explicitly as follows:

\begin{quotation}In a series of independent trials, an event $E$ has the constant probability $p$.  If, in this series, $E$ occurs at least $r$ times in succession, we say that there is a \textbf{\emph{run}} of $r$ successes.  \emph{What is the probability of having a run of $r$ successes in $n$ trials, where naturally $n>r$?}\end{quotation}

Let us denote by $y_{n}$ the unknown probability of a run of $r$ in $n$ trials.  The classical solution to the problem, (\textsc{Feller} \cite{Feller}, \textsc{Uspensky} \cite{Uspensky}), consists of deducing the following \emph{difference equation} for the \emph{complementary} probability $z_{n}:=1-y_{n}$:\begin{equation}
\fbox{$\dis z_{n+1}-z_{n}+qp^{r}z_{n-r}=0$}
\end{equation} where $q:=1-p$, and then concluding that the \emph{generating function}:\begin{equation}
G(x):=z_{0}+z_{1}x+z_{2}x^{2}+\cdots+z_{n}x^{n}+[\cdots]
\end{equation}is, in fact, the \emph{rational} function \begin{equation}
\fbox{$\dis G(x)=\frac{1-p^{r}x^{r}}{1-x+qp^{r}x^{r+1}}$}
\end{equation}
The coefficient of $x^{n}$ gives the general formula for $z_{n}$.
We will prove:\begin{thm}Let\begin{equation}
\beta_{n,r}:=\sum_{l=0}^{[\frac{n}{r+1}]}(-1)^{l}\binom{n-lr}{l}(qp^{r})^{l}
\end{equation}where $[\alpha]:=$ the greatest integer contained in $\alpha$. Then
\begin{equation}
\fbox{$\dis z_{n}=\beta_{n,r}-p^{r}\beta_{n-r,r}$}
\end{equation}
\end{thm}

This explicit formula for $z_{n}$, which one would expect to be at least as famous as the problem, itself, is amazingly hard to find in the literature.  Although \textsc{Todhunter} \cite{Todhunter} details the solutions of \textsc{De Moivre}, \textsc{Condorcet}, and \textsc{Laplace}, none of them gives the formula, although they all give versions of the generating function.

Nor is it to be found in the classical text of \textsc{Markoff} \cite{Markoff}.

When we looked at the modern texts of \textsc{Chung} \cite{Chung}, \textsc{Feller} \cite{Feller},  \textsc{Gnedenko} \cite{Gnedenko}, \textsc{Hoel} \cite{Hoel}, \textsc{Parzen} \cite{Parzen}, \textsc{Ross} \cite{Ross}, \textsc{Tucker} \cite{Tucker}, and \textsc{Uspensky} \cite{Uspensky}, we were able to find only a \emph{statement}, and that \textbf{\emph{without proof}}, of \textbf{Theorem 1}, only in \textsc{Uspensky} \cite{Uspensky}.  The most detailed presentation of the theory of runs on the internet is to be found in \textsc{Weisstein} \cite{Weisstein}, but the explicit formula is not even mentioned there! 

\textsc{Uspensky}, (\cite{Uspensky}, page 79), states that the formula can be found ``. . . according to the known rules." Following \textsc{Feller} (\cite{Feller}, pp. 275-276), if we write \begin{equation}
G(x)=\frac{U(x)}{V(x)}
\end{equation}
then the formula for $z_{n}$ is:\begin{equation}
z_{n}=\frac{\rho_{1}}{x_{1}^{n+1}}+\frac{\rho_{2}}{x_{2}^{n+1}}+\cdots+\frac{\rho_{r+1}}{x_{r+1}^{n+1}}
\end{equation}
where $x_{1}, x_{2}, \cdots, x_{r+1}$ are the $r+1$ distinct roots of $V(x)=0$ and \begin{equation}
\rho_{k}=\frac{-U(x_{k})}{V'(x_{k})}.
\end{equation}
Unfortunately, the equation $V(x)=0$ cannot, in general, be solved explicitly for its $r+1$ roots, and so ``. . . the known rules" are useless in this case.

We therefore offer the following simple derivation of the formula for $z_{n}$, based on the binomial theorem.

\section{Proof of Theorem 1}

It suffices to prove the formula for $\beta_{n,r}$ since the formula for $z_{n}$ is an immediate consequence of it. Now,\begin{align*}
    \sum_{n=0}^{\infty}\beta_{n,r}x^{n}&\equiv
  \frac{1}{1-x+qp^{r}x^{r+1}} \\&= \frac{1}{1-x(1-qp^{r}x^{r})}  \\
    &=\sum_{k=0}^{\infty}\{x(1-qp^{r}x^{r})\}^{k}\\
    &=\sum_{k=0}^{\infty}x^{k}\sum_{l=0}^{k}\binom{k}{l}(-qp^{r}x^{r})^{l}\\
    &=\sum_{k=0}^{\infty}\sum_{l=0}^{k}(-1)^{l}\binom{k}{l}(qp^{r})^{l}x^{rl+k}.\\
    \end{align*}

We must now determine the coefficient, $\beta_{n,r}$, of $x^{n}$ in this series.  Thus, for fixed $r$ and $n$, we must find all pairs of integers $(l,k)$, with $0\leqslant l\leqslant k$ which satisfy the equation: \begin{equation}
rl+k=n,
\end{equation}since each such pair contributes the summand\begin{equation}(-1)^{l}\binom{k}{l}(qp^{r})^{l}=(-1)^{l}\binom{n-lr}{l}(qp^{r})^{l}\end{equation}to the final value of $\beta_{n,r}$.
By inspection we note that\begin{align*}
\label{}
    n&=r\cdot 0+n   \\
    &=r\cdot 1+(n-r)\\ 
    &=r\cdot 2+(n-2r)\\
     &=r\cdot 3+(n-3r)\\
      &=\cdots \cdots\\
       &=(n-r)+r\cdot 1.\\
\end{align*}and these equations correspond to the pairs $$
(l,k)=(0,n), \ (1,n-r), \ (2,n-2r), \cdots, \ (n-r,r),$$respectively.  We do not include $(l,k)=(n,0)$ since the corresponding summand has the value $0$.

By (9) and the final eqution in the list above, the largest value of $l$ occurs when $l=k$ and thus satisfies the equation:\begin{equation*}
lr+r=n
\end{equation*}and we conclude that\begin{equation}
l=\left[\frac{n}{r+1}\right].
\end{equation}Therefore, $l$ takes on the values $0, \ 1, \ 2,\ \cdots, \ \left[\frac{n}{r+1}\right]$, and the proof is complete.


\end{document}